\documentclass[11pt,leqno]{article}
\usepackage{amsmath, amsthm, amsfonts, amssymb}
\usepackage{mathrsfs}
\newfam\msbfam
\font\tenmsb=msbm10    \textfont\msbfam=\tenmsb \font\sevenmsb=msbm7
\scriptfont\msbfam=\sevenmsb \font\fivemsb=msbm5
\scriptscriptfont\msbfam=\fivemsb

\newfam\bigfam
\font\tenbig=msbm10 scaled \magstep2   \textfont\bigfam=\tenbig
\font\sevenbig=msbm7 scaled \magstep2 \scriptfont\bigfam=\sevenbig
\font\fivebig=msbm5 scaled \magstep2
\scriptscriptfont\bigfam=\fivebig

\def\dint{\displaystyle\int}

\def\dsum{\displaystyle\sum}
\def\dfrac{\displaystyle\frac}

\begin{document}


\begin{center} {\Large\bf  $p$-adic singular integral and their commutator in generalized Morrey space}

\vspace{0.3cm}

\noindent{\large\bf Huixia Mo\footnote[1]{
Corresponding author.\\
E-mail addresses: huixiamo@bupt.edu.cn.\\}, Zhe Han, Liu Yang}\\
School of Science, Beijing University
of Posts and Telecommunications, Beijing 100876,  China

\end{center}

\parskip 12pt

\noindent{\bf Abstract}\;\;{\footnotesize
 For a prime number $p,$ let $\mathbb{Q}_p$ be the field of $p$-adic numbers. In this paper, we  established the boundedness of  a class
of $p$-adic singular integral operators on  the $p$-adic generalized  Morrey  spaces.
The corresponding boundedness for the commutators generalized by  the $p$-adic  singular integral operators and $p$-adic Lipschitz functions or $p$-adic generalized Campanato functions is also considered.}
 \vspace{0.3cm}


\section{\bf Introduction}

Let $p$ be a prime number and $x\in\mathbb{Q}.$  Then the non-Archimedean $p$-adic normal $|x|_p$  is defined as follows: if $x=0,$ $|0|_p=0;$ if $x\neq 0$ is an arbitrary rational number with the unique representation $x=p^{\gamma}\frac{m}{n},$ where $m, n$ are not divisible by $p,$ $\gamma=\gamma(x)\in \mathbb{Z},$ then $|x|_p=p^{-\gamma}.$ This normal satisfies $|xy|_p=|x|_p|y|_p,$ $|x+y|_p\leq\max\{|x|_p, |y|_p\}$ and
$|x|_p=0$ if and only if $x=0.$ Moreover, when  $|x|_p\neq|y|_p,$ we have $|x+y|_p=\max\{|x|_p, |y|_p\}.$
 Let $\mathbb{Q}_p$ be the field of $p$-adic numbers, which is  defined as the  completion of the field of rational numbers $\mathbb{Q}$ with respect to the non-Archimedean $p$-adic normal $|\cdot|_p.$
For $\gamma\in\mathbb{Z},$ we denote the ball $B_{\gamma}(a)$   with  center at $a\in\mathbb{Q}_p^n$ and radius $p^\gamma$
and its boundary $S_{\gamma}(a)$  by
 $$\begin{array}{cl}B_{\gamma}(a)=\{x\in\mathbb{Q}_p^n: |x-a|_p\leq p^{\gamma}\},\;\;\;
 S_{\gamma}(a)=\{x\in\mathbb{Q}_p^n: |x-a|_p=p^{\gamma}\},\end{array}$$
respectively.  It is easy to see that
  $$B_\gamma(a)=\bigcup\limits_{k\leq \gamma}S_{k}(a).$$

For $n\in\mathbb{N},$ the space $\mathbb{Q}_p^n=\mathbb{Q}_p\times\cdots\times\mathbb{Q}_p$ consists of all points ${x}=(x_1,\cdots, x_n)$ where $x_i\in \mathbb{Q}_p,$ $i=1,\dots,n,$ $n\geq1.$ The $p$-adic norm of $\mathbb{Q}_p^n$ is defined by
             $$|{x}|_p=\max\limits_{1\leq i\leq n}{|x_i|_p},\;\; {x}\in\mathbb{Q}_p^n.$$
Thus, it is easy to see that $|{x}|_p$ is a non-Archimedean norm on $\mathbb{Q}_p^n.$
              The  balls $B_\gamma({a})$ and the sphere $S_\gamma({a})$ in $\mathbb{Q}_p^n,$ $\gamma\in\mathbb{Z}$ are defined similar to the case $n=1.$

     Since $\mathbb{Q}_p^n$ is a locally compact commutative group under addition, thus from the standard analysis  there  exists the Haar measure $dx$ on the additive group  $\mathbb{Q}_p^n$ normalized by   $\int_{B_0}dx=|B_0|_{H}=1,$  where $|E|_H$ denotes the Haar measure of a measurable set $E\subset\mathbb{Q}_p^n.$
     Then by a simple calculation the Haar measures of any balls and spheres can be
obtained. From the integral theory, it is easy to see that  $|B_{\gamma}({a})|_H=p^{n\gamma}$ and $|S_{\gamma}({a})|_H=p^{n\gamma}(1-p^{-n})$ for any ${a}\in\mathbb{Q}_p^n.$
 For a more complete introduction to the $p$-adic analysis, one can refer to \cite{VVZ,S,R,CEKMM,T,H2,K,WF} and the references therein.

The $p$-adic numbers have been applied in string theory, turbulence theory, statistical mechanics,
quantum mechanics, and so forth(see \cite{VVZ,AK,ABK} for detail). In the past few years, there is an increasing interest in the study of harmonic analysis on $p$-adic field (see \cite{T,H2,K,WF} for detail).

Let $\Omega\in L^{\infty}(\mathbb{Q}_p^n),$  $\Omega(p^jx)=\Omega(x)$ for all
$j\in \mathbb{Z}$ and $\int_{|x|_p=1}\Omega(x)dx=0.$ Then the $p$-adic singular integral
operator defined by Taibleson \cite{T} is as follows
$$T_k(f)(x)=\dint_{|y|_p>p^k}f(x-y)\dfrac{\Omega(y)} {|y|_p^n}dz.$$
 And the $p$-adic singular integral operator $T$ is defined as the limit of $T_k$ when $k$ goes to
 $-\infty.$

Moreover, let$\overrightarrow{b}=(b_1,b_2,...,b_m),$  where $b_i\in L_{loc}{(\mathbb{Q}_p^n)}$ for $1\leq i\leq m.$  Then the higher commutator generated by $\vec{b}$ and $T_k$  can be defined by
$$T_k^{\vec{b}}f(x)=\dint_{|y|_p>p^k}\prod\limits_{i=1}^{m}(b_i(x)-b_i(x-y))f(x-y)\dfrac{\Omega(y)} {|y|_p^n}dz.
$$
 And the commutator generated by  $\overrightarrow{b}=(b_1,b_2,...,b_m)$ and $p$-adic singular integral operator $T$ is defined as the limit of $T_k^{\vec{b}},$ when $k$ goes to  $-\infty.$

Under some conditions, the authors in \cite{T,PT}, obtained that $T_k$
 were of type $(q, q), 1 < q < \infty, $ and of weak
type $(1,1)$ on local fields. In \cite{WLF},  Wu etal. established the
boundedness of $T_k$  on p-adic central Morrey spaces. Furtherly,  the $\lambda$-central BMO
estimates for commutators of these singular integral operators on $p$-adic central Morrey spaces were
obtained in \cite{WLF}. Moreover, in $p$-adic linear space $\mathbb{Q}_p^n$ , Volosivets \cite {V} gave the sufficient
conditions for the maximal function and Riesz potential in $p$-adic generalized Morrey spaces.  Mo etal.\cite{MWM} established the boundedness of the commutators  generated by the $p$-adic Riesz potential  and $p$-adic generalized Campanato functions  in $p$-adic generalized Morrey spaces.

Motivated by the works of \cite{WLF,V,MWM},
we are going to consider  the boundedness of $T_k$
 on the $p$-adic generalized  Morrey type spaces, as well as the boundedness of the commutators
 generated by  $L_k$ and $p$-adic generalized Campanato functions.

Throughout this paper, the letter $C $  will be used to denote various constants and the various uses of
the letter do not, however, denote the same constant. And, $A\lesssim B$ means that $A\leq CB,$ with some
 positive constant $C$ independent of appropriate quantities.

\section{\bf Some notations and  lemmas}

\textbf{Definition 2.1}\cite{V}~Let  $1\leq q<\infty,$ and let $\omega(x)$ be a non-negative measurable function
in $\mathbb{Q}_p^n .$
 A function $f\in L^q_{loc}(\mathbb{Q}_p^n )$ is said to belong to the generalized Morrey space  $GM_{q,\omega}(\mathbb{Q}_p^n ),$ if
$$\|f\|_{GM_{q,\omega}}=\sup\limits_{a\in\mathbb{Q}^n_p,\gamma\in\mathbb{Z}}\dfrac{1}{\omega(B_{\gamma}(a))}\biggl(\dfrac{1}{|B_{\gamma}(a)|_H}\dint_{B_{\gamma}(a)}
|f(y)|^{q}dy\biggr)^{1/q}<\infty,$$
where $\omega(B_{\gamma}(a))=\int_{B_{\gamma}(a)}\omega(x)dx.$

Let $\lambda\in \mathbb{R}.$ If $\omega(B_{\gamma}(a))=|B_{\gamma}(a)|^{\lambda},$ then  $GM_{q,\omega}(\mathbb{Q}_p^n )$ is the classical Morrey spaces $M_{q,\lambda}(\mathbb{Q}_p^n ).$
About the  generalized Morrey space, see \cite{N}, and the classical Morrey spaces, see {\cite{M}}, etc.

Moreover, let $\lambda\in \mathbb{R}$ and $1\leq q<\infty.$  The $p$-adic central Morrey space $CM_{q,\lambda}(\mathbb{Q}^n_p)$ (see \cite{WF}) is defined by
 $$\|f\|_{CM_{q,\lambda}}=\sup\limits_{\gamma\in\mathbb{Z}}\biggl(\dfrac{1}{|B_\gamma(0)|_H^{1+\lambda q}}\dint_{B_\gamma(0)}
|f(y)|^{q}dy\biggr)^{1/q}<\infty.$$

\textbf{Definition 2.2}\cite{CD}~\;\;Let $0<\beta<1$ , then the $p$-adic Lipschitz space $\Lambda_{\beta}(\mathbb{Q}^{n}_p)$ is defined
the set of all functions $f: \mathbb{Q}_p^n\mapsto\mathbb{C}$ such that

$$\|f\|_{\Lambda_{\beta}(\mathbb{Q}_p^n)}=\sup\limits_{x,h\in\mathbb{Q}_p^n, h\neq 0}
\frac{|f(x+h)-f(x)|}{|h|^{\beta}}<\infty.$$

\textbf{Definition 2.3}\cite{V}~Let $B$ be a ball in $\mathbb{Q}_p^n ,$ $1\leq q<\infty.$ And let $\omega(x)$ be a non-negative measurable  function
in $\mathbb{Q}_p^n .$  A function $f\in L^q_{loc}(\mathbb{Q}_p^n )$ is said to belong to the
generalized Campanato space  $GC_{q,\omega}(\mathbb{Q}_p^n ),$ if
$$\|f\|_{GC_{q,\omega}}=\sup\limits_{a\in\mathbb{Q}^n_p,\gamma\in\mathbb{Z}}\dfrac{1}{\omega(B_{\gamma}(a))}\biggl(\dfrac{1}{|B_{\gamma}(a)|_H}
\dint_{B_{\gamma}(a)}|f(y)-f_{B_{\gamma}(a)}|^{q}dy\biggr)^{1/q}<\infty,$$
where $f_{B_{\gamma}(a)}=\frac{1}{|B_{\gamma}(a)|_{H}}\int_{B}f(x)dx$ and $\omega(B_{\gamma}(a))=\int_{B_{\gamma}(a)}\omega(x)dx.$

The  classical Campanato spaces can be seen in \cite{C}, \cite{P} and etc.
The important particular case of $GC_{q,\omega}(\mathbb{Q}_p^n )$ is CBMO$_{q,\lambda}(\mathbb{Q}_p^n ),$ where for $1<q<\infty$ and $ 0<\lambda<1/n.$ And the central BMO space
CBMO$_{q,\lambda}(\mathbb{Q}_p^n )$ is defined by
$$\|f\|_{CBMO^{q,\lambda}(\mathbb{Q}_p^n )}=\sup\limits_{\gamma\in\mathbb{Z}}\dfrac{1}{|B_\gamma(0)|_H^\lambda}
\biggl(\dfrac{1}{|B_\gamma(0)|_H}\dint_{B_\gamma(0)}|f(y)-f_{B_\gamma(0)}|^{q}dy\biggr)^{1/q}<\infty.\eqno{(2.1)}$$

\textbf{Lemma 2.1}\cite{MWM}\;\;Let $1\leq q<\infty,$ amd let $\omega$ be a non-negative measurable   function.
Suppose that $b\in GC_{q,\omega}(\mathbb{Q}_p^n ),$ then
$$|b_{B_k(a)}-b_{B_{j}(a)}|\leq\|b\|_{GC_{q,\omega}}|j-k|\max\{\omega(B_k(a)),\omega(B_j(a))\},$$
for $j,k\in\mathbb{Z}$ and any fixed $a\in\mathbb{Q}_p^n.$\\

Thus, for $j>k,$ from Lemma 2.1, it deduce that
$$\biggl(\dint_{B_{j}(a)}|b(y)-b_{B_{k}(a)}|^{q}dy\biggr)^{1/q}
\leq(j+1-k)|B_j(a)|_{H}^{1/q}\omega(B_j(a))\|b\|_{GC_{q,\omega}}.\eqno{(2.2)}$$

\textbf{Lemma 2.2} \cite{T}\;\; Let $\Omega\in L^{\infty}(\mathbb{Q}_p^n),$  $\Omega(p^jkx)=\Omega(x)$ for all
$j\in \mathbb{Z}$ and $\dint_{|x|_p=1}\Omega(x)dx=0.$ If
$$\sup\limits_{|y|_p=1}\sum\limits_{j=1}^{\infty}\dint_{|x|_p=1}|\Omega(x+p^jy)-\Omega(x)|dx<\infty,$$
then for $1 < p <\infty $, there is a constant $C > 0$
 such that
 $$\|T_k(f)\|_{L^p(\mathbb{Q}_p^n)}\leq C\|f\|_{L^p(\mathbb{Q}_p^n)},$$
for $k\in \mathbb{Z}$, where $C$ is independent of $f$ and $k \in {\mathbb Z}.$

Furthermore, $T(f) = \lim\limits_{k\rightarrow -\infty}T_k(f)$
 exists in the $L^p$ norm and
$$\|T(f)\|_{L^p(\mathbb{Q}_p^n)}\leq C\|f\|_{L^p(\mathbb{Q}_p^n)}.$$

Moreover, on the $p-$adic field,  the Riesz potential $I_\alpha^{p}$ is defined by
$$I^\alpha_{p}f(x)=\dfrac{1}{\Gamma_n(\alpha)}\dint_{\mathbb{Q}_p^n}\dfrac{f(y)}{|x-y|_p^{n-\alpha}}dy,$$
where $\Gamma_n(\alpha)=(1-p^{\alpha-n})/(1-p^{-\alpha}),$ $\alpha\in \mathbb{C}$ and $\alpha\neq0.$

\textbf{Lemma 2.3}\cite{MWM}\;\;Let $\alpha$ be a complex number with $0< {\mbox Re}\alpha<n,$ and let  $1<r<\infty,$ $1<q<n/ {\mbox Re}\alpha,$ $0<1/r=1/q- {\mbox Re}\alpha/n.$
 Suppose that both $\omega$ and $\nu$ are non-negative measurable  functions, such that
$$\dsum_{j={\gamma}}^{\infty}p^{j {\mbox Re}\alpha}\dfrac{\nu(B_j(a))}{\omega(B_\gamma(a))}=C<\infty,$$
for any $a\in\mathbb{Q}^n_p$ and  $\gamma\in\mathbb{Z}.$ Then the Riesz potential $I^{\alpha}_p$ is bounded from $GM_{q,\nu}$ to $GM_{r,\omega}.$\\

\section{\bf Main results}
  In this section, let us state the main results of the paper.

\textbf{Theorem 3.1}\;\;
\;\; Let $1<q<\infty,$ and let  $\Omega(p^jx)=\Omega(x)$ for all
$j\in \mathbb{Z},$ $\dint_{|x|_p=1}\Omega(x)dx=0,$ and
$$\sup\limits_{|y|_p=1}\sum\limits_{j=1}^{\infty}\dint_{|x|_p=1}|\Omega(x+p^jy)-\Omega(x)|dx<\infty.$$
 Suppose that both $\omega$ and $\nu$ are non-negative measurable  functions, such that
$$\dsum_{j={\gamma}}^{\infty}\nu(B_j(a))/\omega(B_\gamma(a))=C<\infty,\eqno{(3.1)}$$
for any $\gamma\in\mathbb{Z}$ and $a\in\mathbb{Q}^n_p.$ Then the singular integral operators $T_k$ are bounded from $GM_{q,\nu}$ to $GM_{q,\omega}$ for all
$k\in \mathbb{Z}.$ Moreover, $T(f)=\lim\limits_{k\rightarrow-\infty}T_k(f)$  exists in  $GM_{q,\omega}$ and the operator $T$ is
 bounded from $GM_{q,\nu}$ to $GM_{q,\omega}.$  \\

\textbf{Corollary 3.1}\;\; Let $1<q<\infty,$  $\lambda<0,$  $\Omega\in L^{\infty}(\mathbb{Q}_p^n),$  $\Omega(p^jx)=\Omega(x)$ for all
$j\in \mathbb{Z},$ $\dint_{|x|_p=1}\Omega(x)dx=0,$ and
$$\sup\limits_{|y|_p=1}\sum\limits_{j=1}^{\infty}\dint_{|x|_p=1}|\Omega(x+p^jy)-\Omega(x)|dx<\infty.$$
Then the operators $T_k$ and $T$ are bounded on the space $CM_{q,\lambda}$ for all
$k\in \mathbb{Z}.$

In fact, for $\lambda<0.$ Taking $\omega(B)=\nu(B)=|B|_H^{\lambda}$ in Theorem 3.1, we can obtain the Corollary 3.1.   If  the Morrey space $M^{q,\lambda}(\mathbb{Q}_p^n )$ is replaced by the central Morrey space $CM^{ q,\lambda}(\mathbb{Q}_p^n )$ in Corollary 3.1, then the conclusion is  that of Theorem 4.1 in \cite{WLF}.

\textbf{Theorem 3.2}\;\; Let $\Omega\in L^{\infty}(\mathbb{Q}_p^n),$  $\Omega(p^jx)=\Omega(x)$ for all
$j\in \mathbb{Z},$ $\dint_{|x|_p=1}\Omega(x)dx=0,$ and
$$\sup\limits_{|y|_p=1}\sum\limits_{j=1}^{\infty}\dint_{|x|_p=1}|\Omega(x+p^jy)-\Omega(x)|dx<\infty.$$
Let  $0<\beta_i<1$ for $i=1,2,\dots,m,$ such that $0<\beta=\dsum_{i=1}^{m}\beta_i<n.$ And, let $1<r<\infty,$ $1<q<n/\beta$ such that  $1/r=1/q-\beta/n.$
 Suppose that $b_i\in {\Lambda_{\beta_i}},i=1,2,\dots,m,$ and  both $\omega$ and $\nu$ are non-negative measurable  functions, such that
$$\dsum_{j={\gamma}}^{\infty}p^{j\beta}\nu(B_j(a))/\omega(B_\gamma(a))=C<\infty,\eqno{(3.2)}$$
for any $\gamma\in\mathbb{Z}$ and $a\in\mathbb{Q}^n_p.$ Then the commutators $T_k^{\vec{b}}$ are  bounded from $GM_{q,\nu}$ to $GM_{r,\omega},$ for all $k\in\mathbb{Z}.$ Moreover,  $T^{\vec{b}}(f)=\lim\limits_{k\rightarrow-\infty}T_k^{\vec{b}}(f)$ exists in the space of $GM_{q,\omega},$ and the commutator $T^{\vec{b}}$ is bounded from $GM_{q,\nu}$ to $GM_{q,\omega}.$

\textbf{Theorem 3.3}\;\;
\;\; Let $\Omega\in L^{\infty}(\mathbb{Q}_p^n),$  $\Omega(p^jx)=\Omega(x)$ for all
$j\in \mathbb{Z},$ $\dint_{|x|_p=1}\Omega(x)dx=0,$ and
$$\sup\limits_{|y|_p=1}\sum\limits_{j=1}^{\infty}\dint_{|x|_p=1}|\Omega(x+p^jy)-\Omega(x)|dx<\infty.$$
Let $1< q,r,q_1,\dots q_m<\infty,$ such that $1/r=1/q+1/q_{1}+1/q_{2}+\dots+1/q_{m}.$
 Suppose that $\omega,$ $\nu$ and $\nu_i$ $ (i=1,2,\dots,m)$ are non-negative measurable functions.
If $b_i\in GC_{q_i,\nu_i}(\mathbb{Q}_p^n )$ for $i=1,2,\dots,m,$
and  the functions  $\omega,$ $\nu$ and $\nu_i(i=1,2,\dots,m)$  satisfy the following conditions\\
(i)\;\;$\prod\limits_{i=1}^m\nu_i(B_\gamma(a))\nu(B_\gamma(a))/\omega(B_\gamma(a))=C<\infty,$\\
(ii)\;\;$\dsum_{j=\gamma+1}^{\infty}\prod\limits_{i=1}^m\nu_i(B_j(a))(j+1-\gamma)^m\nu(B_j(a))/\omega(B_\gamma(a))=C<\infty,$\\
for any $\gamma\in\mathbb{Z}$ and $a\in\mathbb{Q}^n_p.$
Then the  commutators  $T_k^{\vec{b}}$  are bounded from $GM_{q,\nu}$ to $GM_{r,\omega},$ for all
$k\in \mathbb{Z}.$ Moreover, the commutator $T^{\vec{b}}=\lim\limits_{k\rightarrow-\infty}T_k^{\vec{b}}$ exists in the space of $GM_{q,\omega},$  and
$T^{\vec{b}}$ is bounded from $GM_{q,\nu}$ to $GM_{q,\omega}.$

\textbf{Corollary 3.2}\;\;
\;\; Let $\Omega\in L^{\infty}(\mathbb{Q}_p^n),$  $\Omega(p^jx)=\Omega(x)$ for all
$j\in \mathbb{Z},$ $\dint_{|x|_p=1}\Omega(x)dx=0,$ and
$$\sup\limits_{|y|_p=1}\sum\limits_{j=1}^{\infty}\dint_{|x|_p=1}|\Omega(x+p^jy)-\Omega(x)|dx<\infty.$$
Let $1< q,r,q_1,\dots q_m<\infty,$ such that $1/r=1/q+1/q_{1}+1/q_{2}+\dots+1/q_{m}.$ Let $0\leq \lambda_1,\dots,\lambda_m<1/n,$
 $\lambda<-\sum\limits_{i=1}^{m}\lambda_i$ and $\tilde{\lambda}=\sum\limits_{i=1}^{m}\lambda_i+\lambda.$
If $b_i\in BMO_{q_i,\lambda_i}(\mathbb{Q}_p^n ),$
then the commutators  $T_k^{\vec{b}}$ and $T^{\vec{b}}$ are bounded from $M_{q,\lambda}$ to $M_{r,\tilde{\lambda}}.$

Moreover, let $1<r,q,q_1<\infty,$ such that $1/r=1/q+1/q_{1}.$ Let $0\leq \lambda_1<1/n,$
 $\lambda<-\lambda_1$ and $\tilde{\lambda}=\lambda_1+\lambda.$
 If $b\in CBMO_{q_1,\lambda_1}(\mathbb{Q}_p^n ),$ then from Corollary 4.1, it follows that the commutator $[T_k, b]$ and $[T, b]$ both are bounded from
 $CM_{q,\lambda}$ to $CM_{r,\tilde{\lambda}}.$  This conclusion is  that of  Theorem 4.2 in \cite{WLF}.

\section{\bf Proof of Theorem 3.1-3.3}

Let us give the proof of  Theorem 3.1, firstly.\\

 For any fixed $\gamma\in \mathbb{Z}$ and  $a\in\mathbb{Q}^n_p,$ it is easy to see that
$$\begin{array}{cl}
 &\dfrac{1}{\omega(B_{\gamma}(a))}\biggl(\dfrac{1}{|B_{\gamma}(a)|_H}\dint_{B_{\gamma}(a)}|T_k(f)(x)|^qdx\biggr)^{1/q}\\
 \leq&\dfrac{1}{\omega(B_{\gamma}(a))}\biggl(\dfrac{1}{|B_{\gamma}(a)|_H}\dint_{B_{\gamma}(a)}|T_k(f)(f\chi_{B_\gamma(a)})(x)|^qdx\biggr)^{1/q}\\
 &+\dfrac{1}{\omega(B_{\gamma}(a))}\biggl(\dfrac{1}{|B_{\gamma}(a)|_H}\dint_{B_{\gamma}(a)}|T_k(f\chi_{B^{c}_\gamma(a)})(x)|^qdx\biggr)^{1/q}\\
:=&I+II,
\end{array}\eqno{(4.1)}$$
where $B^{c}_\gamma(a)$ is the complement to $B_\gamma(a)$ in  $\mathbb{Q}^n_p.$

Using Lemma 2.2 and  (3.1), it follows that
$$\begin{array}{cl}
 I=&\dfrac{1}{\omega(B_{\gamma}(a))}\biggl(\dfrac{1}{|B_{\gamma}(a)|_H}\dint_{B_{\gamma}(a)}|T_k(f\chi_{B_\gamma(a)})(x)|^qdx\biggr)^{1/q}\\
 \leq&\dfrac{1}{\omega(B_{\gamma}(a))}\dfrac{1}{|B_{\gamma}(a)|_H^{1/q}}\biggl(\dint_{B_{\gamma}(a)}|f(x)|^qdx\biggr)^{1/q}\\
=&\dfrac{\nu(B_{\gamma}(a))}{\omega(B_{\gamma}(a))}\dfrac{1}{\nu(B_{\gamma}(a))}\biggl(\dfrac{1}{|B_{\gamma}(a)|_H}\dint_{B_{\gamma}(a)}|f(x)|^qdx\biggr)^{1/q}\\
 \lesssim &\|f\|_{GM_{q,\nu}}.
\end{array}\eqno{(4.2)}$$

For $II,$ let us estimate $|T_k(f\chi_{B^{c}_\gamma(a)})(x)|,$ firstly.

Since $x\in B_{\gamma}(a)$ and $\Omega\in L^{\infty}(\mathbb{Q}_p^n),$ then we have

$$\begin{array}{cl}
 |T_k(f\chi_{B^{c}_\gamma(a)})(x)|
 &=\biggl|\dint_{|y|_p>p^k}(f\chi_{B^{c}_\gamma(a)})(x-y)\dfrac{\Omega(y)} {|y|_p^n}dy \biggr|\\
 &=\biggl|\dint_{|x-z|_p>p^k}(f\chi_{B^{c}_\gamma(a)})(z)\dfrac{\Omega(x-z)} {|x-z|_p^n}dz \biggr|\\
   &\lesssim\dint_{B^{c}_{\gamma}(a)}\dfrac{|f(z)|} {|x-z|_p^n}dz \\
&\lesssim\dsum_{j=\gamma+1}^{\infty}\dint_{S_j(a)}p^{-jn}|f(y)|dy\\
&\leq\dsum_{j=\gamma+1}^{\infty}p^{-jn}\biggl(\dint_{B_j(a)}|f(y)|^qdy\biggr)^{1/q}|B_j(a)|_H^{1-1/q}\\
 &= \|f\|_{GM_{q,\nu}}\dsum\limits_{j=\gamma+1}^{\infty}\nu(B_j(a)).
\end{array}\eqno{(4.3)}$$

Thus, from  (3.1) and (4.3), it follows that

$$\begin{array}{cl}
 II=&\dfrac{1}{\omega(B_\gamma(a))}\biggl(\dfrac{1}{|B_\gamma(a)|_H}\dint_{B_\gamma(a)}|T_k(f\chi_{B^{c}_\gamma(a)})(x)|^qdx\biggr)^{1/q}\\
 \lesssim&\|f\|_{GM_{q,\nu}}\sum\limits_{j=\gamma+1}^{\infty}\nu(B_j(a))/\omega(B_\gamma(a))\\
   \lesssim &\|f\|_{GM_{q,\nu}}
\end{array}\eqno{(4.4)}$$

Combining the estimates of (4.1), (4.2) and (4.4),  we have

$$\begin{array}{cl}
 &\dfrac{1}{\omega(B_{\gamma}(a))}\biggl(\dfrac{1}{|B_{\gamma}(a)|_H}\dint_{B_{\gamma}(a)}|T_k(f)(x)|^qdx\biggr)^{1/q} \lesssim \|f\|_{GM_{q,\nu}},
\end{array}$$
which means that
$T_k$ is bounded from $GM_{q,\nu}$ to $GM_{q,\omega}.$

Moreover, from Lemma 2.2 and the definition of  $GM_{q,\omega}(\mathbb{Q}_p^n ),$ it is obvious that $T(f)=\lim\limits_{k\rightarrow-\infty}T_k(f)$  exists in  $GM_{q,\omega}$ and the operator $T$ is
 bounded from $GM_{q,\nu}$ to $GM_{q,\omega}.$ \\

Proof of Theorem 3.2.\\

 For any $x\in\mathbb{Q}^n_p$, since $\Omega\in L^{\infty}(\mathbb{Q}_p^n),$ and $b_i\in {\Lambda_{\beta_i}},i=1,2,\dots,m,$ then it is easy to see that
$$\begin{array}{cl}&|T_k^{\vec{b}}f(x)|\\
\leq&\dint_{|y|_p>p^k}\prod\limits_{i=1}^{m}|b_i(x)-b_i(x-y)||f(x-y)|\dfrac{|\Omega(y)|} {|y|_p^n}dy\\
\lesssim&\dint_{\mathbb{Q}_p^n}\dfrac{|f(z)|}{|x-z|_p^{n-\beta}}dz\\
\lesssim&I^\beta_{p}(|f|)(x).\end{array}$$

Thus,  from Lemma 2.2 it is obvious that  the commutators $T_k^{\vec{b}}$ are bounded from $GM_{q,\nu}$ to $GM_{r,\omega},$ for all $k\in\mathbb{Z}.$

Moreover, from the definition of  $GM_{q,\omega}(\mathbb{Q}_p^n ),$ it is obvious that $T^{\vec{b}}(f)=\lim\limits_{k\rightarrow-\infty}T_k^{\vec{b}}(f)$    exists in the space of $GM_{q,\omega},$ and the commutator
 $T^{\vec{b}}$ is bounded from $GM_{q,\nu}$ to $GM_{q,\omega}.$  \\

Proof of Theorem 3.3\\
 Without loss of generality, we need only to show that the conclusion holds for
$m=2.$

For any fixed $\gamma\in \mathbb{Z}$ and $a\in \mathbb{Q}^n,$ we write $f^0=f\chi_{B_\gamma(a)}$ and
$f^\infty=f\chi_{B^{c}_\gamma(a)},$ then
$$\begin{array}{cl}
 &\dfrac{1}{\omega(B_\gamma(a))}\biggl(\dfrac{1}{|B_\gamma(a)|_H}\dint_{B_\gamma(a)}|T_k^{(b_1,b_2)}(f)(x)|^rdx\biggr)^{1/r}\\
 \leq&\dfrac{1}{\omega(B_\gamma(a))}\biggl(\dfrac{1}{|B_\gamma(a)|_H}\dint_{B_\gamma(a)}
 |(b_1(x)-(b_1)_{B_\gamma(a)})(b_2(x)-(b_2)_{B_\gamma(a)})T_k(f^0)(x)|^rdx\biggr)^{1/r}\\
 &+\dfrac{1}{\omega(B_\gamma(a))}\biggl(\dfrac{1}{|B_\gamma(a)|_H}\dint_{B_\gamma(a)}
 |(b_1(x)-(b_1)_{B_\gamma(a)})T_k((b_2-(b_2)_{B_\gamma(a)})f^0)(x)|^rdx\biggr)^{1/r}\\
  &+\dfrac{1}{\omega(B_\gamma(a))}\biggl(\dfrac{1}{|B_\gamma(a)|_H}\dint_{B_\gamma(a)}
 |(b_2(x)-(b_2)_{B_\gamma(a)})T_k((b_1-(b_1)_{B_\gamma(a)})f^0)(x)|^rdx\biggr)^{1/r}\\
 &+\dfrac{1}{\omega(B_\gamma(a))}\biggl(\dfrac{1}{|B_\gamma(a)|_H}\dint_{B_\gamma(a)}
 |T_k((b_1-(b_1)_{B_\gamma(a)})(b_2-(b_2)_{B_\gamma(a)})f^0)(x)|^rdx\biggr)^{1/r}\\
  &+\dfrac{1}{\omega(B_\gamma(a))}\biggl(\dfrac{1}{|B_\gamma(a)|_H}\dint_{B_\gamma(a)}
 |(b_1(x)-(b_1)_{B_\gamma(a)})(b_2(x)-(b_2)_{B_\gamma(a)})T_k(f^{\infty})(x)|^rdx\biggr)^{1/r}\\
  &+\dfrac{1}{\omega(B_\gamma(a))}\biggl(\dfrac{1}{|B_\gamma(a)|_H}\dint_{B_\gamma(a)}
 |(b_1(x)-(b_1)_{B_\gamma(a)})T_k((b_2-(b_2)_{B_\gamma(a)})f^{\infty})(x)|^rdx\biggr)^{1/r}\\
 &+\dfrac{1}{\omega(B_\gamma(a))}\biggl(\dfrac{1}{|B_\gamma(a)|_H}\dint_{B_\gamma(a)}
 |(b_2(x)-(b_2)_{B_\gamma(a)})T_k((b_1-(b_1)_{B_\gamma(a)})f^{\infty})(x)|^rdx\biggr)^{1/r}\\
  &+\dfrac{1}{\omega(B_\gamma(a))}\biggl(\dfrac{1}{|B_\gamma(a)|_H}\dint_{B_\gamma(a)}
 |T_k((b_1-(b_1)_{B_\gamma(a)})(b_2-(b_2)_{B_\gamma(a)})f^{\infty})(x)|^rdx\biggr)^{1/r}\\
=:&E_1+E_2+E_3+E_4+E_5+E_6+E_7+E_8.
\end{array}\eqno{(4.5)}$$

In the following, we will estimate every part, respectively.

Since  $1/r=1/q+1/q_{1}+1/q_{2},$ then, from H\"{o}lder's inequality, Lemma 2.2 and (i), it follows that
$$\begin{array}{cl}
E_1=&\dfrac{1}{\omega(B_\gamma(a))}\biggl(\dfrac{1}{|B_\gamma(a)|_H}\dint_{B_\gamma(a)}
 |(b_1(x)-(b_1)_{B_\gamma(a)})(b_2(x)-(b_2)_{B_\gamma(a)})T_k(f^0)(x)|^rdx\biggr)^{1/r}\\
 \leq&\dfrac{1}{\omega(B_\gamma(a))|B_\gamma(a)|_H^{1/r}}\prod\limits_{i=1}^2\biggl(\dint_{B_\gamma(a)}
 |b_i(x)-(b_i)_{B_\gamma(a)}|^{q_i}dx\biggr)^{1/q_{i}}
 \biggl(\dint_{B_\gamma(a)}|T_k(f^0)(x)|^{q}dx\biggr)^{1/{q}}\\
 \leq&\dfrac{\nu_1(B_\gamma(a))\nu_2(B_\gamma(a))}{\omega(B_\gamma(a))|B_\gamma(a)|_H^{1/q}}\prod\limits_{i=1}^2\|b_i\|_{GC_{q_i,\nu_i}}
 \biggl(\dint_{B_\gamma(a)}|f(x)|^{q}dx\biggr)^{1/q}\\
 \leq&\dfrac{\nu(B_\gamma(a))\nu_1(B_\gamma(a))\nu_2(B_\gamma(a))}{\omega(B_\gamma(a))}\prod\limits_{i=1}^2\|b_i\|_{GC_{q_i,\nu_i}}
 \|f\|_{GM_{q,\nu}}\\
\lesssim&\prod\limits_{i=1}^{2}\|b_i\|_{GC_{q_i,\nu_i}}\|f\|_{GM_{q,\nu}}.
\end{array}$$

Let $1/\bar{q}=1/q+1/q_{2},$  then $1/r=1/q_{1}+1/\bar{q}.$ Thus, from H\"{o}lder's inequality, Lemma 2.2 and (i), we obtain
$$\begin{array}{cl}
E_2&=
\dfrac{1}{\omega(B_\gamma(a))}\biggl(\dfrac{1}{|B_\gamma(a)|_H}\dint_{B_\gamma(a)}
 |(b_1(x)-(b_1)_{B_\gamma(a)})T_k((b_2-(b_2)_{B_\gamma(a)})f^0)(x)|^rdx\biggr)^{1/r}\\
\leq&\dfrac{1}{\omega(B_\gamma(a))|B_\gamma(a)|_H^{1/r}}\biggl(\dint_{B_\gamma(a)}
 |b_1(x)-(b_1)_{B_\gamma(a)}|^{q_1}dx\biggr)^{1/q_{1}}
 \biggl(\dint_{B_\gamma(a)}|T_k((b_2-(b_2)_{B_\gamma(a)})f^0)(x)|^{\bar{q}}dx\biggr)^{1/{\bar{q}}}\\
 \leq&\dfrac{1}{\omega(B_\gamma(a))|B_\gamma(a)|_H^{1/r}}\biggl(\dint_{B_\gamma(a)}
 |b_1(x)-(b_1)_{B_\gamma(a)}|^{q_1}dx\biggr)^{1/q_{1}}
 \biggl(\dint_{B_\gamma(a)}|(b_2(x)-(b_2)_{B_\gamma(a)})f(x)|^{\bar{q}}dx\biggr)^{1/{\bar{q}}}\\
 \leq&\dfrac{1}{\omega(B_\gamma(a))}\prod\limits_{i=1}^2\biggl(\dint_{B_\gamma(a)}
 |b_i(x)-(b_i)_{B_\gamma(a)}|^{q_i}dx\biggr)^{1/q_{i}} \biggl(\dint_{B_\gamma(a)}|f(x)|^{q}dx\biggr)^{1/{q}}\\
  \leq&\dfrac{\nu(B_\gamma(a))\nu_1(B_\gamma(a))\nu_2(B_\gamma(a))}{\omega(B_\gamma(a))}
  \prod\limits_{i=1}^2\|b_i\|_{GC_{q_i,\nu_i}} \|f\|_{GM_{q,\nu}}\\
 \lesssim&\prod\limits_{i=1}^{2}\|b_i\|_{GC_{q_i,\nu_i}}\|f\|_{GM_{q,\nu}}.
\end{array}$$

Similarly,
$$\begin{array}{cl}
E_3 \lesssim\prod\limits_{i=1}^2\|b_i\|_{GC_{q_i,\nu_i}}
 \|f\|_{GM_{q,\nu}}.\\
\end{array}$$

For $E_4,$ from Lemma 2.2, H\"{o}lder's inequality and  (i), we obtain
$$\begin{array}{cl}
E_4&=
\dfrac{1}{\omega(B_\gamma(a))}\biggl(\dfrac{1}{|B_\gamma(a)|_H}\dint_{B_\gamma(a)}
|T_k(b_1-(b_1)_{B_\gamma(a)})(b_2-(b_2)_{B_\gamma(a)})f^0)(x)|^rdx\biggr)^{1/r}\\
\leq&\dfrac{1}{\omega(B_\gamma(a))|B_\gamma(a)|_H^{1/r}}
 \biggl(\dint_{B_\gamma(a)}|(b_1(x)-(b_1)_{B_\gamma(a)})(b_2(x)-(b_2)_{B_\gamma(a)})f(x)|^{r}dx\biggr)^{1/{r}}\\
  \leq&\dfrac{1}{\omega(B_\gamma(a))|B_\gamma(a)|_H^{1/r}}\prod\limits_{i=1}^2\biggl(\dint_{B_\gamma(a)}
 |b_i(x)-(b_i)_{B_\gamma(a)}|^{q_i}dx\biggr)^{1/q_{i}} \biggl(\dint_{B_\gamma(a)}|f(x)|^{q}dx\biggr)^{1/{q}}\\
  \leq&\dfrac{\nu(B_\gamma(a))\nu_1(B_\gamma(a))\nu_2(B_\gamma(a))}{\omega(B_\gamma(a))}\prod\limits_{i=1}^2\|b_i\|_{GC_{q_i,\nu_i}}
 \|f\|_{GM_{q,\nu}}\\
  \lesssim&\prod\limits_{i=1}^{2}\|b_i\|_{GC_{q_i,\nu_i}}\|f\|_{GM_{q,\nu}}.
\end{array}$$

To estimate $E_5,$ we need to consider $|T_k(f^{\infty})(x)|,$ firstly. In fact from (4.3), it is easy to see that
$$
 |T_k(f^{\infty})(x)|\lesssim\|f\|_{GM_{q,\nu}}\dsum\limits_{j=\gamma+1}^{\infty}\nu(B_j(a)).\eqno{(4.6)}$$
 
Therefore, from H\"{o}lder's inequality, (4.6) and (ii), we obtain
$$\begin{array}{cl}
E_5=&\dfrac{1}{\omega(B_\gamma(a))}\biggl(\dfrac{1}{|B_\gamma(a)|_H}\dint_{B_\gamma(a)}
 |(b_1(x)-(b_1)_{B_\gamma(a)})(b_2(x)-(b_2)_{B_\gamma(a)})T_k(f^{\infty})(x)|^rdx\biggr)^{1/r}\\
  \leq&\dfrac{1}{\omega(B_\gamma(a))|B_\gamma(a)|_H^{1/r}}\prod\limits_{i=1}^2\biggl(\dint_{B_\gamma(a)}
 |b_i(x)-(b_i)_{B_\gamma(a)}|^{q_i}dx\biggr)^{1/q_{i}} \biggl(\dint_{B_\gamma(a)}|T_k(f^{\infty})(x)f(x)|^{q}dx\biggr)^{1/{q}}\\
\lesssim&\dsum_{j=\gamma+1}^{\infty}\dfrac{\nu(B_j(a))\nu_1(B_\gamma(a))\nu_2(B_\gamma(a))}{\omega(B_\gamma(a))}
 \prod\limits_{i=1}^2\|b_i\|_{GC_{q_i,\nu_i}}
 \|f\|_{GM_{q,\nu}}\\
\lesssim&\prod\limits_{i=1}^{2}\|b_i\|_{GC_{q_i,\nu_i}}\|f\|_{GM_{q,\nu}}.
\end{array}$$

It is similar to the estimate of (4.3), for $x\in B_{\gamma}(a),$ by  $\Omega\in L^{\infty}(\mathbb{Q}_p^n)$ and (2.2), we can deduce that
$$\begin{array}{cl}
 &|T_k(b_2-(b_2)_{B_{\gamma}(a)})f^{\infty})(x)|\\
 &=\biggl|\dint_{|y|_p>p^k}(b_2(x-y)-(b_2)_{B_{\gamma}(a)(a)})f\chi_{B^{c}_\gamma(a)}(x-y)\dfrac{\Omega(y)} {|y|_p^n}dy \biggr|\\
 &\leq\dint_{B^{c}_{\gamma}}|b_2(z)-(b_2)_{B_{\gamma}(a)}||f(z)|\dfrac{|\Omega(x-z)|} {|x-z|_p^n}dz \\
  &\lesssim\dint_{B^{c}_{\gamma}}\dfrac{|b_2(z)-(b_2)_{B_{\gamma}(a)}||f(z)|} {|x-z|_p^n}dz \\
&\lesssim\dsum_{j=\gamma+1}^{\infty}\dint_{S_j(a)}p^{-jn}|b_2(z)-(b_2)_{B_{\gamma}(a)}||f(y)|dy\\
&=\dsum_{j=\gamma+1}^{\infty}p^{-jn}|B_j(a)|_H^{1-1/q-1/q_2}\biggl(\dint_{S_j(a)}|f(y)|^{q}dy\biggr)^{1/{q}}
\biggl(\dint_{S_j(a)(a)}|b_2(y)-(b_2)_{B_{\gamma}(a)}|^{q_2}dy\biggr)^{1/q_2}\\
&\leq\|f\|_{GM_{q,\nu}}\dsum_{j=\gamma+1}^{\infty}p^{-jn}
|B_j(a)|_H^{1-1/q_2}\nu(B_j(a))\biggl(\dint_{B_j(a)}|b_2(y)-(b_2)_{B_{\gamma}(a)}|^{q_2}dy\biggr)^{1/q_2}\\
&\lesssim\|b_2\|_{GC_{q_2,\nu_2}}\|f\|_{GM_{q,\nu}}\dsum_{j=\gamma+1}^{\infty}(j+1-\gamma)\nu(B_j(a))\nu_2(B_j(a)).
\end{array}\eqno{(4.7)}$$

Let $1/\bar{q}=1/q+1/q_{2},$  then $1/r=1/q_{1}+1/\bar{q}.$ Thus, from H\"{o}lder's inequality, (4.7) and  (ii), it follows that
$$\begin{array}{cl}
E_6=&\dfrac{1}{\omega(B_\gamma(a))}\biggl(\dfrac{1}{|B_\gamma(a)|_H}\dint_{B_\gamma(a)}
 |(b_1(x)-(b_1)_{B_\gamma(a)})T_k((b_2-(b_2)_{B_\gamma(a)})f^{\infty})(x)|^rdx\biggr)^{1/r}\\
  \leq&\dfrac{1}{\omega(B_\gamma(a))|B_\gamma(a)|_H^{1/r}}\biggl(\dint_{B_\gamma(a)}
 |b_1(x)-(b_1)_{B_\gamma(a)}|^{q_1}dx\biggr)^{1/q_{1}} \biggl(\dint_{B_\gamma(a)}|T_k((b_2-(b_2)_{B_\gamma(a)})f^{\infty})(x)|^{\bar{q}}dx\biggr)^{1/{\bar{q}}}\\
 \leq&\prod\limits_{i=1}^2\|b_i\|_{GC_{q_i,\nu_i}}
 \|f\|_{GM_{q,\nu}}\dfrac{1}{\omega(B_\gamma(a))}\dsum_{j=\gamma+1}^{\infty}(j+1-\gamma)\nu(B_j(a))\nu_2(B_j(a))\nu_1(B_\gamma(a))\\
\lesssim&\prod\limits_{i=1}^{2}\|b_i\|_{GC_{q_i,\nu_i}}\|f\|_{GM_{q,\nu}}.
\end{array}$$

It's analogues to the estimate of $E_6,$ we obtain

 $$\begin{array}{cl}
E_7\lesssim\prod\limits_{i=1}^{2}\|b_i\|_{GC_{q_i,\nu_i}}\|f\|_{GM_{q,\nu}}.
\end{array}$$

Moreover, since $\Omega\in L^{\infty}(\mathbb{Q}_p^n),$  then by (2.2), it is easy to see that
$$\begin{array}{cl}
 &|T_k((b_1-(b_1)_{B_{\gamma}(a)})(b_2-(b_2)_{B_{\gamma}(a)})f^{\infty})(x)|\\
 &=\biggl|\dint_{|x-z|_p>p^k}(b_1(z)-(b_1)_{B_{\gamma}(a)})(b_2(z)-(b_2)_{B_{\gamma}(a)})f\chi_{B^{c}_\gamma(a)}(z)\dfrac{\Omega(x-z)} {|x-z|_p^n}dz \biggr|\\
 &\leq\dint_{B^{c}_{\gamma}}|b_1(z)-(b_1)_{B_{\gamma}(a)}||b_2(z)-(b_2)_{B_{\gamma}(a)}||f(z)|\dfrac{|\Omega(x-z)|} {|x-z|_p^n}dz \\
 &\lesssim\dsum_{j=\gamma+1}^{\infty}\dint_{S_j(a)}p^{-jn}|b_1(z)-(b_1)_{B_{\gamma}(a)}||b_2(z)-(b_2)_{B_{\gamma}(a)}||f(y)|dy\\
&=\dsum_{j=\gamma+1}^{\infty}p^{-jn}|B_j(a)|_H^{1-1/q-1/q_1-1/q_2}\biggl(\dint_{S_j(a)}|f(y)|^{q}dy\biggr)^{1/{q}}
\biggl(\dint_{S_j(a)}|b_1(y)-(b_1)_{B_{\gamma}(a)}|^{q_1}dy\biggr)^{1/q_1}\\
&\quad\times\biggl(\dint_{S_j(a)}|b_2(y)-(b_2)_{B_{\gamma}(a)}|^{q_2}dy\biggr)^{1/q_2}\\
&\lesssim\prod\limits_{i=1}^{2}\|b_i\|_{GC_{q_i,\nu_i}}\|f\|_{GM_{q,\nu}}\dsum_{j=\gamma+1}^{\infty}(j+1-\gamma)^2\nu(B_j(a))\nu_1(B_j(a))\nu_2(B_j(a)).
\end{array}\eqno{(4.8)}$$

Therefore, from (4.8) and (ii), we get that
$$\begin{array}{cl}
E_8&=\dfrac{1}{\omega(B_\gamma(a))}\biggl(\dfrac{1}{B_\gamma(a)|_H}\dint_{B}
 |T_k((b_1-(b_1)_{B_\gamma(a)})(b_2-(b_2)_{B_\gamma(a)})f^{\infty})(x)|^rdx\biggr)^{1/r}\\
 &\leq\prod\limits_{i=1}^{2}\|b_i\|_{GC_{q_i,\nu_i}}\|f\|_{GM_{q,\nu}}\dfrac{1}{\omega(B_\gamma(a))}\dsum_{j=\gamma+1}^{\infty}(j+1-\gamma)^2
 \nu(B_j(a))\nu_1(B_j(a))\nu_2(B_j(a))\\
 &\lesssim\prod\limits_{i=1}^{2}\|b_i\|_{GC_{q_i,\nu_i}}\|f\|_{GM_{q,\nu}}.
\end{array}$$

Combining (4.5) and the estimates of $E_1,E_2,\dots, E_8,$  we have
$$\begin{array}{cl}
 &\dfrac{1}{\omega(B_\gamma(a))}\biggl(\dfrac{1}{|B_\gamma(a)|_H}\dint_{B_\gamma(a)}|T_k^{(b_1,b_2)}(f)(x)|^rdx\biggr)^{1/r}
 \leq\prod\limits_{i=1}^{2}\|b_i\|_{GC_{q_i,\nu_i}}\|f\|_{GM_{q,\nu}},
\end{array}$$
which means that the  commutator $T_k^{(b_1,b_2)}$ is bounded from $GM_{q,\nu}$ to $GM_{r,\omega}.$

Moreover, from Lemma 2.2 and the definition of  $GM_{q,\omega}(\mathbb{Q}_p^n ),$ it is obvious that $T^{\vec{b}}(f)=\lim\limits_{k\rightarrow-\infty}T_k^{\vec{b}}(f)$  exists in the space of $GM_{q,\omega},$ and the commutator
 $T^{\vec{b}}$ is bounded from $GM_{q,\nu}$ to $GM_{q,\omega}.$ \\

 Therefore, the proof of Theorem 3.3 is complete.

\noindent{\bf Acknowledgements}\;\;{This work is supported by the National Natural Science Foundation of China (11601035)}

\end{document}